\documentclass[twoside,11pt,leqno]{article}

\usepackage{amsmath, amssymb}
\usepackage{enumerate}
\usepackage[all]{xy}
\CompileMatrices

\let\bsbf=\bf

\def\iff{if and only if }
\def\qed{\hfill$\Box$}

\newtheorem{theorem}{Theorem}[section]
\newtheorem{lemma}[theorem]{Lemma}

\newcommand{\rnf}{\renewcommand{\thefootnote}{\arabic{footnote}}}

\def\cpx{\mathbb{C}}

\def\SU{\operatorname{SU}}
\def\SO{\operatorname{SO}}

\def\R{{\mathbb R}}

\def\Z{\mathbb{Z}}

\title{\vspace{-3cm}
\ \hspace{-2.70in}
{\small {\bf Morfismos},
Vol. 13, No. 2, 2009,
pp.
        55--\pageref{ultimapagina} 
}
\\
\vspace{2.6cm} 
$K$-theory of $S^7/Q_8$ and a counterexample to a result of P.M. Akhmet'ev
\footnote{\bsbf Invited article.}
}

\author{\rnf   Peter S.~Landweber
}

\date{}

\addtocounter{page}{54}

\begin{document}
\maketitle

\begin{abstract} \noindent
A simple counterexample is presented to a proposition which is used in the arguments given by P.~M.~Akhmet'ev in his work on the Hopf invariant and Kervaire invariant.  The counterexample makes use of the $K$-theory of the quotient of the $7$-sphere by the quaternion group of order $8.$
\end{abstract}

\noindent 

\section{The counterexample}

In Petr Akhmet'ev's paper \cite{Ak1} on the Hopf invariant, Proposition 41 (Proposition 37 in the English translation) in \S 4 is a key step in the proof of the Main Theorem, which is given at the end of \S 2.  The same Proposition is essential to the success of the proof of the Main Theorem in Akhmet'ev's paper \cite{Ak2} on the Kervaire invariant, which is given in \S 5 of that paper; more precisely, Lemma 22 in \S 5 (in the English translation) is used in the proof of the Main Theorem, and this lemma depends on the Proposition in the paper \cite{Ak1}.  In both papers, the use of this Proposition enters toward the end of the proofs of the Main Theorems, and does not concern the parts of the arguments viewed as the main steps by Akhmet'ev.  

\smallskip
We shall give the statement of Proposition 41 of \cite{Ak1} and then produce a simple counterexample, which makes use of the $K$-theory of the quotient space $S^7/Q_8$.

\smallskip
In order to state Proposition 41 of \cite{Ak1}, we begin by recalling the character table 
for the quaternion group $Q_8$ of order $8$ from 
\cite[\S 13]{At1}, where $x_i, x_j$, and $x_k$ are 1-dimensional complex representations and $y$ is a $2$-dimensional complex representation arising from the inclusion 
$$
Q_8 \subset S^3 = \SU(2).
$$
Here, $Q_8 = \{\pm 1, \pm i, \pm j, \pm k \} \subset \mathbb{H}.$
We view $Q_8$ as acting on the unit sphere in 
$\mathbb{H} \oplus \mathbb{H}$ 
(which we identify with 
$\mathbb{C}^2 \oplus \mathbb{C}^2$)
by left quaternion multiplication.
To each complex representation $V$ of $Q_8$ we associate a complex vector bundle on $S^7/Q_8$ with total space the orbit space
$(S^7 \times V)/Q_8.$  Let $\xi_i, \xi_j, \xi_k,$ and $\eta$ denote the complex vector bundles on $S^7/Q_8$ obtained from 
$x_i, x_j, x_k,$ and $y$, respectively.  For a complex vector bundle 
$\zeta$, we denote its underlying oriented real vector bundle by 
$\zeta_{\R}.$  For example, $\eta_{\R}$ is an oriented real $4$-plane bundle on $S^7/Q_8.$  Notice that $S^7/Q_8$ is orientable, since $Q_8$ is contained in the connected group $\SU(2).$  

\vspace{-6mm}
  
\centerline{\begin{picture}(160,176)(-25,-72)
\put(40,79){\small (Conjugacy classes)}
\put(-76,4){\small (Irreducible}
\put(-86,-8){\small representations)}
\put(0,0){\begin{tabular}{c||c|c|c|c|c}
\rule{0pt}{15pt} & $\,\,1\,\,$ & $-1$ & $\pm i$ & $\pm j$ & $\pm k$ \\[.5ex]
\hline\hline
\rule{0pt}{15pt} $1$ & $1$ & $1$ & $1$ &$1$ & $1$ \\[.5ex]
\hline
\rule{0pt}{15pt} $x_i$ & $1$ & $1$ & $1$ & $-1$ & $-1$ \\[.5ex]
\hline
\rule{0pt}{15pt} $x_j$ & $1$ & $1$ & $-1$ & $1$ & $-1$ \\[.5ex]
\hline
\rule{0pt}{15pt} $x_k$ & $1$ & $1$ & $-1$ & $-1$ & $1$ \\[.5ex]
\hline
\rule{0pt}{15pt} $y$ & $2$ & $-2$ & $0$ & $0$ & $0$ \\[.5ex]
\end{tabular}}
\end{picture}}

Proposition 41 in \cite[\S 4]{Ak1} is the following assertion:
\textit{Let $f: K^7 \to S^7/Q_8$ be a continuous map of a smooth closed oriented $7$-manifold to $S^7/Q_8$, and assume that the stable normal bundle to $K^7$ has the form $f^*(2k\, \eta_{\R})$ with $k$ odd.  Then the degree of $f$ is even.} 

\smallskip
We shall show that this is not true, as a consequence of the following two results.

\begin{theorem}
	The tangent bundle $\tau$ to $S^7/Q_8$ satisfies

\vspace{-3mm}\ 
	$$
	\tau \oplus \varepsilon^1_{\R} \cong \eta_{\R} \oplus \eta_{\R}
	= 2\, \eta_{\R}.
	$$
\end{theorem}

\vspace{1mm}Here we denote by $\varepsilon^k_{\R}$ the trivial real $k$-plane bundle, with similar notation for trivial complex bundles.

\begin{theorem}
	In $K(S^7/Q_8), \, 2 - \eta$ has order $8$.
\end{theorem}

Hence $8\eta = \varepsilon^{16}_{\mathbb{C}}$ in complex K-theory, so 
$8\, \eta_{\R} = \varepsilon^{32}_{\R}$ in real K-theory.  It follows that the stable normal bundle to $S^7/Q_8$ is $6\, \eta_{\R}$, and so has the form $2k\, \eta_{\R}$ with $k=3.$  So the identity map of 
$S^7/Q_8$ is a counterexample to Proposition 41 in \cite{Ak1}.

\smallskip
Theorem 1.1 is an instance of Theorem 3.1 in Szczarba's paper \cite{Sz}, which concerns the tangent bundle of a quotient $S^n/G$ of a sphere for the case in which a finite group $G$ acts freely on the sphere; in turn, this result is an instance of \cite[Theorem 1.1]{Sz}.  In the next section, we give a direct proof of Theorem 1.1, and also present several proofs of Theorem~1.2.

I am indebted to Petr Akhmet'ev for discussions over the last two months leading to this understanding of \S 4 of his paper [1], and to Jes\'us Gonz\'alez for help with the character tables.

\section{The $K$-theory of $S^7/Q_8$}

\noindent \emph{Proof of Theorem 1.1.}
Let $W$ be a real representation of a finite group $G$ such that $G$ acts freely on the unit sphere $S(W).$ 
Note first that, as $G$-vector bundles on $S(W)$ we have for the total space $T(S(W))$ of the tangent bundle to $S(W)$

\vspace{-3mm}\ 
$$
T\big(S(W)\big) \oplus N\big(S(W)\big)\, \cong \,S(W) \times W,
$$
where $N(S(W))$ is the total space of the normal bundle to $S(W)$ in $W.$  Using outward unit normals to produce an equivariant and nonvanishing section of the normal bundle, we see that 
$N(S(W)) \cong S(W) \times \R$ as $G$-vector bundles, where $G$ acts trivially on the factor $\R.$  Passing to vector bundles on $S(W)/G$ (with no $G$-action), we obtain Szczarba's result 
\cite[Theorem 3.1]{Sz}:

\vspace{-3mm}\ 
$$
T\big(S(W)/G\big) \,\oplus\, S(W)/G \times \R \,\cong \,\big(S(W)\times W\big)/G.
$$
Theorem 1.1 follows by taking $G$ to be $Q_8$ and $V$ to be  
$\mathbb{H} \oplus \mathbb{H} = \cpx^2 \oplus \cpx^2$ with the action being given by left quaternion multiplication. \qed

\medskip
\noindent
{\it First proof of Theorem~1.2 (in a weaker form).}
In positive even dimen\-sions the integral cohomology of $S^7/Q_8$ consists of $2$-torsion, since $S^7/Q_8$ and the classifying space $BQ_8$ have the same cohomology through dimension $6.$  Therefore the Atiyah-Hirzebruch spectral sequence implies that 
$\widetilde{K}(S^7/Q_8)$ is a finite 2-group, so that $2 - \eta$ has order $2^i$ for some integer $i \geq 0.$  The possibility that $i \leq 2$ is ruled out by characteristic class considerations (see the following lemma), so $i \geq 3.$  (By Theorem 1.1 the stable normal bundle to $S^7/Q_8$ is $(2^i-2)\, \eta_{\R} = 
2(2^{i-1}-1)\, \eta_{\R}$ with $i \geq 3,$ which has the form $2k\, \eta_{\R}$ with $k$ odd, and so suffices to establish the counterexample.) \qed

\begin{lemma} A generator of $H^4(S^7/Q_8; \Z) \cong \Z/8$ can be chosen to be the second Chern class $c_2(\eta).$  Hence, 
$c_2(2k\eta)=0$  \iff  $k \equiv 0 \pmod{4}.$  
\end{lemma}

\vspace{-2.5mm}
\noindent
{\it Proof.}
By tom Dieck \cite[pp.~188--189]{tD}, this group is cyclic of order~$8$ 
and a generator can be taken to be the Euler class of the $\SO(4)$-bundle $\eta_{\R}.$  Since this Euler class is equal to $c_2(\eta)$ and $c_1(\eta) = 0$ because $\eta$ is an $\SU(2)$-bundle, 
$c_2(2k\eta) = 2k c_2(\eta)$ and the result follows.\qed

\medskip
\noindent \emph{Second proof of Theorem 1.2.\!}  
The paper \cite{Fu} of Kens\^{o} Fujii contains (in Theorem 1.2) 
a computation of 
$\widetilde{K}(S^{4n+3}/Q_8)$ as an abelian group for all 
$n \geq 0,$ including the result that in 
$\widetilde{K}(S^7/Q_8)$ the element $2-\eta$ generates a cyclic summand of order $8.$  (The multiplicative structure is also determined.) \qed

\medskip \noindent \emph{Third proof of Theorem 1.2.} 
As indicated by Atiyah \cite{At2} in his review of \cite{Ma}, which is available on MathSciNet, if $V$ is a complex representation of a finite group $G$ such that $G$ acts freely on the unit sphere $S(V),$ then the complex $K$-theory of the quotient manifold $S(V)/G$ can be computed by the formula

\vspace{-2.5mm}\ 
$$
K\big(S(V)/G\big) \,\cong \,R(G)/\big(\lambda_{-1}(V)\big),
$$
where $R(G)$ is the representation ring of $G$ and 
$\lambda_{-1}(V) = \sum (-1)^q\lambda^q(V)$ with $\lambda^q$ the $q$th exterior power.  (This is immediate from the Gysin sequence for 
$V \to \mathrm{pt}$ in equivariant $K$-theory.)

\smallskip
In the case of interest here, we have 
$V = \mathbb{C}^2 \oplus \mathbb{C}^2$ corresponding to the representation $2y = y \oplus y$ of $Q_8$.  Since $y$ arises from the inclusion of $Q_8$ into $SU(2),$ $\lambda_2(y) = 1$ and so 
$\lambda_{-1}(2y) = 2-y.$  Following \cite[\S 13]{At1}, we introduce the notation

\vspace{-6mm}\ 
$$
\alpha = 1-x_i,\quad \beta = 1-x_j,\quad \gamma = 3-x_i -x_j -x_k,\quad \delta = 2-y,
$$

\vspace{-1mm}\noindent 
a basis over $\mathbb{Z}$ for the augmentation ideal $I(Q_8).$  Hence, 
$\lambda_{-1}(y) = \delta$ and $\lambda_{-1}(2y) = \delta ^2.$  Note that $2-\eta$ corresponds to $2-y = \delta$ under the isomorphism
$K(S^7/Q_8) \cong R(Q_8)/(\delta^2).$

\smallskip
Our task is to determine $R(Q_8)/(\delta^2);$ for practice we shall also determine $R(Q_8)/(\delta).$  It is convenient to extend the character table for $Q_8$ to include the four virtual representations just introduced 
and to note the effect of multiplication by $\delta$ on them: 

\vspace{-3.5mm}\ 
$$
\alpha\delta=2\alpha,\quad \beta\delta=2\beta,\quad 
\gamma\delta=2\gamma,\quad \delta^2 = 4\delta - \gamma
$$

\vspace{-1.5mm}\noindent
(the last of these involves a slight correction to \cite[\S 13]{At1}).

\centerline{\begin{picture}(160,147)(-25,-59)
\put(40,68){\small (Conjugacy classes)}
\put(-70,4){\small (Virtual}
\put(-86,-8){\small representations)}
\put(0,0){\begin{tabular}{c||c|c|c|c|c}
\rule{0pt}{15pt} & $\,\,1\,\,$ & $-1$ & $\pm i$ & $\pm j$ & $\pm k$ \\[.5ex]
\hline\hline
\rule{0pt}{15pt} $\alpha$ & $0$ & $0$ & $0$ & $2$ & $2$ \\[.5ex]
\hline
\rule{0pt}{15pt} $\beta$ & $0$ & $0$ & $2$ & $0$ & $2$ \\[.5ex]
\hline
\rule{0pt}{15pt} $\gamma$ & $0$ & $0$ & $4$ & $4$ & $4$ \\[.5ex]
\hline
\rule{0pt}{15pt} $\delta$ & $0$ & $4$ & $2$ & $2$ & $2$ \\[.5ex]
\end{tabular}}
\end{picture}}

The image of $\delta$ acting by multiplication on $R(Q_8)$ is easily seen to be generated by $2\alpha,\, 2\beta,\, \gamma,$ and $\delta.$  Hence $\widetilde{K}(S^3/Q_8)$ is isomorphic to 
$\mathbb{Z}/2 \oplus \mathbb{Z}/2$ with generators given by $\alpha$ and $\beta.$  

\smallskip
Finally, the image of $\delta^2$ acting by multiplication on $R(Q_8)$ is generated by the virtual representations
$4\alpha,\, 4\beta,\, 2\gamma, \textrm \, 4\delta - \gamma,$
so also contains $8\delta.$  It is more convenient to take as generators for $(\delta ^2)$ the virtual representations

\vspace{-6mm}\ 
$$
4\alpha, \quad 4\beta, \quad 8\delta, \quad \gamma - 4\delta,
$$
from which it follows easily that

\vspace{-2mm}\ 
$$ 
\widetilde{K}(S^7/Q_8) \,\cong \,\mathbb{Z}/4 \oplus
\mathbb{Z}/4 \oplus \mathbb{Z}/8
$$
with generators represented by $\alpha,\, \beta,$ and $\delta,$ respectively. (Note that $\gamma \equiv 4\delta$ in this quotient.)  
\qed

\vspace{2mm}
\hfill\
{\footnotesize
\parbox{5.15cm}{ $\rule{2mm}{0mm}$ \\ 
$\rule{2mm}{0mm}$ \\ $\rule{2mm}{0mm}$ \\Peter S.~Landweber \\
{\it Department of Mathematics}\\
Rutgers University\\
Piscataway, NJ 08854\\ 
{\sf landwebe@math.rutgers.edu}}} {\hfill}

\label{ultimapagina}
\end{document}